\documentclass{article}
\usepackage[american]{babel}
\usepackage{amsfonts,amsmath,amssymb,epsf,epsfig}
\usepackage{xypic}
\xyoption{all}
\input{xypic}
\newdir{ >}{{}*!/-8pt/\dir{>}}

\def\R{\mathbb{R}}
\def\N{\mathbb{N}}
\def\C{\mathbb{C}}

\def\Z{\mathbb{Z}}

\def\K{\mathbb{K}} 
    
\def\fg{{\frak g}}

\def\fm{{\frak m}}

\def\pr{\noindent $\bf{Proof.}$\quad}     
\def\fin{\hfill$\square$\\}           
\newtheorem{theo}{Theorem}

\newtheorem{prop}{Proposition}

\begin{document}
\date{}
\title{Cohomology and deformations of the infinite dimensional filiform Lie 
algebra $\fm_2$}

\author{Alice Fialowski\\
        E\"otv\"os University, Budapest\\
	fialowsk@cs.elte.hu\\
\and    Friedrich Wagemann \\
        Universit\'e de Nantes\\
	wagemann@math.univ-nantes.fr}
   
\maketitle


\begin{abstract}
Denote $\fm_2$ the infinite dimensional $\N$-graded Lie algebra defined by 
the basis $e_i$ for $i\geq 1$ and by relations $[e_1,e_i]=e_{i+1}$ for all 
$i\geq 2$, $[e_2,e_j]=e_{j+2}$ for all $j\geq 3$.
We compute in this article the bracket structure on $H^1(\fm_2,\fm_2)$,
$H^2(\fm_2,\fm_2)$ and in relation to this, we establish that there are only
finitely many true deformations of $\fm_2$ in each weight by constructing
them explicitely. It turns out that in weight $0$ one gets as non-trivial 
deformations only one formal non-converging deformation.
\end{abstract}

\section*{Introduction}

Recall the classification of infinite dimensional $\N$-graded Lie algebras 
$\fg=\bigoplus_{i=1}^{\infty}\fg_i$ with one-dimensional homogeneous components
$\fg_i$ and two generators over a field of characteristic zero. A. Fialowski 
showed in \cite{Fia1} that any Lie algebra of this type must be isomorphic to 
$\fm_0$, $\fm_2$ or $L_1$. We call these Lie algebras infinite dimensional 
filiform Lie algebras in analogy with the finite dimensional case where the 
name was coined by M. Vergne in \cite{Ver}. Here $\fm_0$ is given by generators
$e_i$, $i\geq 1$, and relations $[e_1,e_i]=e_{i+1}$ for all $i\geq 2$, $\fm_2$ 
with the same generators by relations $[e_1,e_i]=e_{i+1}$ for all $i\geq 2$, 
$[e_2,e_j]=e_{j+2}$ for all $j\geq 3$, and $L_1$ with the same generators is 
given by the relations $[e_i,e_j]=(j-i)e_{i+j}$ for all $i,j\geq 1$ 
(cf the beginning of section $1$ for our convention about writing relations
for a Lie algebra). $L_1$
appears as the positive part of the Witt algebra given by generators $e_i$ for
$i\in\Z$ with the same relations $[e_i,e_j]=(j-i)e_{i+j}$ for all $i,j\in\Z$.
The result was also obtained later by Shalev and Zelmanov in \cite{SZ}.

The cohomology with trivial coefficients of the Lie algebra $L_1$
was studied in \cite{Gont},
the adjoint cohomology in degrees $1$, $2$ and $3$ has been computed in 
\cite{Fia2} and also all of its non equivalent deformations were given. 
For the Lie algebra $\fm_0$, the cohomology with trivial coefficients 
has been studied in \cite{FialMill}, and the adjoint cohomology in degrees 
$1$ and $2$ in \cite{FiaWag}. The adjoint cohomology in degrees $1$ and $2$
of $\fm_2$ is the object of the present article. The cohomology of $\fm_0$ and 
$\fm_2$ rose interest only recently, and the reason is probably that - as 
happens usually for solvable Lie algebras - the cohomology is huge and 
therefore meaningless. Our point of view is that there still remain interesting
features.   

Indeed, it is true that the first and second adjoint cohomology of $\fm_2$ are 
infinite dimensional, but they are much less impressive than the analoguous 
results for $\fm_0$. We believe that this comes from the much more restrictive
bracket structure for $\fm_2$. Actually, the bracket structure is so rigid
that there is no infinite dimensional filiform Lie algebra ``between''
$\fm_2$ and $L_1$. 
The space $H^1(\fm_2,\fm_2)$ becomes already interesting 
when we split it up into homogeneous components $H^1_l(\fm_2,\fm_2)$ 
of weight $l\in\Z$, this latter space being finite dimensional for each 
$l\in\Z$. The bracket structure on $H^1(\fm_2,\fm_2)$ is studied in section 
$2$.

The space $H^2(\fm_2,\fm_2)$ is discussed in section $3$. This space is here 
finite dimensional in each weight separately. 
Given a generator of $H^2(\fm_2,\fm_2)$, i.e. an infinitesimal deformation,
corresponding to the linear term of a formal deformation, one can try to 
adjust higher order terms in order to satisfy the Jacobi identity in the deformed
Lie algebra up to order $k$. If the Jacobi identity is satisfied to all orders,
we will call it a true (formal) deformation, see Fuchs' book \cite{Fuks} for
details on cohomology and \cite{Fia2} for deformations of Lie algebras.
  
In section $3.2$ we discuss Massey products, in sections $3.3$ -- $3.5$ 
we describe all true deformations in negative weights. Section $3.6$ identifies
the deformations in weight zero.

As obstructions to infinitesimal deformations given by classes in 
$H^2(\fm_2,\fm_2)$ are expressed by Massey powers of these classes in 
$H^3(\fm_2,\fm_2)$, it is the vanishing of these Massey squares, cubes etc
which makes it possible to prolongate an infinitesimal deformation to all 
orders. For $\fm_2$ here, on the one hand the cocycle equations are so rigid
that they select already few cochains to be cocycles, but on the other hand,
there are enough cochains to compensate all Massey powers, leading to formal,
non-converging deformations. The main result reads

\begin{theo}
The true deformations of $\fm_2$ are finitely generated in each weight. 
More precisely, the space of unobstructed cohomology classes is zero
in weight $l\leq -5$, because there are no non-trivial cocycles. It is
in degree $l\geq -4$ of dimension two (but with changing representatives),
but only of dimension one for $l=-1,0,1$, because one cocycle becomes 
a coboundary in these weights.

The infinitesimal deformation in weight $l=0$ can be prolongated to all 
orders and gives a formal non-converging deformation. 
\end{theo}

As a rather astonishing consequence, $\fm_2$ does not deform to $L_1$.

We believe that the discussion of these examples of deformations are 
interesting as they go beyond the usual approach where the condition that    
$H^2({\mathfrak g},{\mathfrak g})$ should be finite dimensional is the 
starting point for the examination of deformations, namely the existence 
of a miniversal deformation \cite{FiaFuc}.

Another attractive point of our study is the fact that here for $\fm_2$ 
the Massey squares, cubes etc. involved can all be compensated and lead to an 
interesting obstruction calculus. Thus the second adjoint cohomology of 
$\fm_2$ may serve as an example on which to study explicitely obstruction 
theory.

After this work has been finished, a preprint of Dimitri 
Millionschikov \cite{Million} appeared, which has much overlap with ours. While he
computes the adjoint cohomology in a more conceptual way using the
Feigin-Fuchs 
spectral sequence, our paper clarifies the bracket structure on $H^1$ and 
the structure of the true deformations. 

\noindent{\bf Acknowledgements:}\quad  
Both authors are grateful to Max Planck Institute in Bonn where this work 
was accomplished. We thank Dimitri Millionschikov for pointing out
an error in some dimensions of the cohomology spaces.
The authors also thank the referee for the careful reading and
suggestions.   

\section{Preliminaries}

This article is about a Lie algebra over a field $\K$ defined below by 
generators and relations; let us specify the ground field $\K$ to be $\R$ or
$\C$, although this does not play a r\^ole in the computations that follow. 
Anyway, we will freely divide by $2$.
  
Recall the $\N$-graded Lie algebra $\fm_2=\bigoplus_{i\geq 1}(\fm_2)_i$; 
all graded components $(\fm_2)_i$ are $1$-dimensional, and we choose a basis
$e_i$ of each of them. The brackets then read: $[e_1,e_i]=e_{i+1}$ for all 
$i\geq 2$, $[e_2,e_j]=e_{j+2}$ for all $j\geq 3$. These relations are always
understood to be the only non-trivial relations (i.e. one has for example also
the relations $[e_3,e_j]=0$ for all $j\geq 3$), except for those which may be 
derived from the given ones by antisymmetry of the bracket.

We will compute in later sections of this paper the Lie algebra cohomology 
spaces $H^1(\fm_2,\fm_2)$ and $H^2(\fm_2,\fm_2)$ of $\fm_2$ with coefficients 
in the adjoint representation. We recommend the book of Dmitry Fuchs 
\cite{Fuks} as a reference on cohomology and deformations, and furthermore
\cite{Fia2} for deformations.
As $\fm_2$ is $\N$-graded, the cochain, cocycle,
coboundary and cohomology spaces are, and thus it makes sense to restrict 
attention to the graded components of {\it weight} $l$ denoted 
$C^*_l(\fm_2,\fm_2)$, $Z^*_l(\fm_2,\fm_2)$, $B^*_l(\fm_2,\fm_2)$ and 
$H^*_l(\fm_2,\fm_2)$ of the spaces of all cochains $C^*(\fm_2,\fm_2)$, cocycles
$Z^*(\fm_2,\fm_2)$, coboundaries $B^*(\fm_2,\fm_2)$ and cohomology classes
$H^*(\fm_2,\fm_2)$.

The cohomology spaces $H^*(\fm_2,\fm_2)$ for $*=1,2$ are interesting from
the following point of view: $H^*(\fm_2,\fm_2)$ carries a graded Lie bracket
$$[,]\,:\,H^p(\fm_2,\fm_2)\otimes H^q(\fm_2,\fm_2)\to H^{p+q-1}(\fm_2,\fm_2),$$
which restricts to a Lie bracket on $H^1(\fm_2,\fm_2)$ which is graded with 
respect to the weight $l$. We will compute this bracket in the next section.

The space $H^2(\fm_2,\fm_2)$ draws its importance from the interpretation of 
being the space of {\it infinitesimal deformations} of the Lie algebra $\fm_2$.
Such an infinitesimal deformation $[\omega]\in H^2(\fm_2,\fm_2)$ is the
term of degree one in the expansion of a deformed bracket with respect to the 
deformation parameter. The question whether the infinitesimal term given by
$[\omega]$ can be prolongated to degree two or even to all higher powers
can be answered by studying the {\it Massey powers} of $[\omega]$. Indeed,
it is a necessary condition for $[\omega]$ to admit a prolongation to degree
two that the {\it Massey square} $[\omega]^2\in H^3(\fm_2,\fm_2)$ is zero, i.e.
if for all $i,j,k\geq 1$
$$\omega(\omega(e_i,e_j),e_k)+{\rm cycl.}\,=\,d\alpha,$$
for some $2$-cochain $\alpha\in C^2(\fm_2,\fm_2)$. In this sense, the Massey 
square is the first obstruction for $[\omega]$ to give a (formal) deformation.
The next obstruction is then the {\it Massey cube}, defined using $\omega$ and
$\alpha$ by
$$\omega(\alpha(e_i,e_j),e_k) + \alpha(\omega(e_i,e_j),e_k)+{\rm cycl.}.$$
In case all obstructions vanish, $[\omega]$ gives rise to a {\it formal 
deformation}. The bracket defined by 
$[,]_t=[,]+t\omega+t^2\alpha+\ldots$ satisfies then the
Jacobi identity up to all orders. But it is not clear whether setting $t=r$ for
some $r\in \R$ defines a Lie bracket $[,]_r$, i.e. it is not clear whether
the formal deformation {\it converges}. If this is the case, we call it a {\it 
true deformation}. A deformation having only a finite number of non-zero
terms is always a true deformation. 

A homogeneous cocycle $\omega$ of weight $l\in\Z$ for the Lie algebras $\fm_0$ 
or $\fm_2$ is given by coefficients $a_{i,j}$ such that 
$\omega(e_i,e_j)=a_{i,j}e_{i+j+l}$. The most important cocycle equation for
$\fm_0$ was (cf \cite{FiaWag}) for $i,j\geq 2$:
$$a_{i+1,j}+a_{i,j+1}\,=\,a_{i,j}.$$
In \cite{FiaWag}, we defined some fundamental solutions to this equation
which we named {\it families}. The $2$-{\it family} has $a_{2,k}=1$ for all 
$k\geq 3$ and $a_{i,j}=0$ for all $i>2$, up to antisymmetry. The 
$3$-{\it family} has $a_{3,k}=1$ for all $k\geq 4$ and $a_{i,j}=0$ for all 
$i>3$, up to antisymmetry. The $a_{2,k}$ coefficients are then easily seen to 
be non-zero starting from $a_{2,5}$, and they grow linearly in $k$. For 
explicit formulae for the $m$-family, we refer to \cite{FiaWag}.        

\section{The space $H^1(\fm_2,\fm_2)$}

We will compute the space $H^1_l(\fm_2,\fm_2)$ of homogeneous cohomology
classes of weight $l\in\Z$ for each fixed $l$. A $1$-cochain 
$\omega\in C^1(\fm_2,\fm_2)$ is called {\it homogeneous of weight $l\in\Z$} in
case $\omega(e_i)\,=\,a_ie_{i+l}$ for each $i\geq 1$. The cocycle identity
reads then for a homogeneous cochain
$$d\omega(e_i,e_j)\,=\,\omega([e_i,e_j])-[e_i,\omega(e_j)]+[e_j,\omega(e_i)]
\,=\,0$$
for all $i,j\geq 1$. We get different sets of equations for $i=1$, $j\geq 2$,
$i=2$, $j\geq 3$, and $i,j\geq 3$.  

\noindent(a)\quad If $i=1$, $j\geq 3$, $j+l\geq 2$: 
$$0\,=\,a_{j+1}-a_j - a_1\delta_{l,0}- a_1\delta_{l,1},$$
if $j\geq 3$, $j+l=0,1$, we get $0\,=\,a_{j+1}$, but there is no equation for 
$j+l\leq -1$.

If $i=1$ and $j=2$, $l\geq 1$:
$$0\,=\,a_3-a_2 + a_1(1-\delta_{l,1}),$$
$0\,=\,a_3-a_2 - a_1$ if $j=2$ and $l=0$, $0\,=\,a_3$ if $j=2$ and $l=-1,-2$,
and no equation if $j+1+l\leq 0$. 
 
\noindent(b)\quad If $i=2$, $j\geq 3$, $j+l\geq 3$:
$$0\,=\,a_{j+2}-a_j-a_2\delta_{l,0}-a_2\delta_{l+1,0},$$
for $j+l=2$, we get $0=a_{-l+4}-a_2\delta_{j,3}$, for $j+l=1$, we get 
$0=a_{-l+3}+a_{-l+1}$, for $j+l=0$, we get $0=a_{-l+2}$, for $j+l=-1$, we get
$0=a_{-l+1}$, and there is no equation for $j+l\leq -2$.

\noindent(c)\quad If $i,j\geq 3$:
$$0\,=\,\delta_{j+l,1}a_j+\delta_{j+l,2}a_j-\delta_{i+l,1}a_i-
\delta_{i+l,2}a_i.$$

Now let us discuss $1$-cocycles in weight $l=0$. For $i=1$ and $j\geq 2$, we 
get by equations (a)
$$0\,=\,a_{j+1}-a_j - a_1,$$
and for $i=2$ and $j\geq 3$ by equations (b)
$$0\,=\,a_{j+2}-a_j - a_2.$$ 
Call $a_1=:a$ and $a_2=:b$, then we get on the one hand $a_3-b=a$, $a_4-a_3=a$,
$a_5-a_4=a$ and so on, and on the other hand $a_5-a_3=b$. Therefore $b=2a$.
In conclusion, we get a one parameter family of cocycles in weight $l=0$.

Now let us discuss $1$-cocycles in weight $l=1$. For $i=1$ and $j\geq 3$, we 
get by equations (a)
$$0\,=\,a_{j+1}-a_j - a_1,$$
while for $j=2$, we get $0=a_3-a_2$. For $i=2$ and $j\geq 3$ by equations (b)
$$0\,=\,a_{j+2}-a_j.$$
We conclude $a_2=a_3$, $a_3=a_5$, $a_1=0$, $a_3=a_4$, and all $a_i$ for 
$i\geq 2$ are then equal. This means that we have one free parameter.

Now let us discuss $1$-cocycles in weight $l\geq 2$. For $i=1$ and $j\geq 3$, 
we get by equations (a)
$$0\,=\,a_{j+1}-a_j,$$
while for $j=2$, we get $0=a_3-a_2+a_1$. For $i=2$ and $j\geq 3$ by equations 
(b)
$$0\,=\,a_{j+2}-a_j.$$
We have $a_4=a_3$ and so on, and $a_1$ and $a_2$ are thus two free parameters.

Now let us discuss $1$-cocycles in weight $l=-1$. For $i=1$ and $j\geq 3$, 
we get by equations (a)
$$0\,=\,a_{j+1}-a_j,$$
while for $j=2$, we get $0=a_3$. For $i=2$ and $j\geq 4$, we get by equations 
(b)
$$0\,=\,a_{j+2}-a_j-a_2,$$
while for $j=3$, we get $0\,=\,a_{5}-a_2$. We have therefore $a_3=0$, 
$a_4=a_3$, $a_5=a_4$, $0=a_6-a_4-a_2$, etc. This gives $a_2=0$, $a_3=0$, 
$a_4=0$, $a_5=0$ and so on. Remark that $a_1$ does not exist, because
$\omega(e_i)=a_ie_{i-1}$. 

Now let us discuss $1$-cocycles in weight $l=-2$. Remark that here $a_1$ and
$a_2$ do not exist. The equations (a), i.e. $i=1$, $j\geq 2$, read
$$0\,=\,a_{j+1}-\left\{\begin{array}{ccc} 0 & {\rm if} & j=2,3 \\
a_j & {\rm if} & j\geq 4 \end{array}\right.$$
The equations (b), i.e. $i=2$, $j\geq 3$, read
$$0\,=\,a_{j+2}+\left\{\begin{array}{ccc} a_3 & {\rm if} & j=3 \\
0 & {\rm if} & j=4 \\
-a_j & {\rm if} & j\geq 5 \end{array}\right.$$
We get thus $a_3=0$, $a_4=0$, $a_5=a_4$, $a_6=0$, and so on. One concludes 
that all coefficients are zero.

Now let us discuss $1$-cocycles in weight $l\leq-3$. Remark that here $a_1$,
$a_2$, up to $a_{-l}$ do not exist. 
The equations (a), i.e. $i=1$, $j\geq 2$, read
$$0\,=\,a_{j+1}-\left\{\begin{array}{ccc} 0 & {\rm if} & j=-l,-l+1 \\
a_j & {\rm if} & j\geq -l+2 \end{array}\right.$$
The equations (b), i.e. $i=2$, $j\geq 3$, read
$$0\,=\,a_{j+2}+\left\{\begin{array}{ccc} 0 & {\rm if} & j= -l-1,-l \\
a_j & {\rm if} & j=-l+1 \\   0 & {\rm if} & j=-l+2 \\
-a_j & {\rm if} & j\geq -l+3 \end{array}\right.$$
One concludes 
that all coefficients are zero.  

Next come the coboundaries. It is clear that $dC^0_l(\fm_2,\fm_2)=0$ for all
weights $l\leq 0$, because coboundaries are brackets with elements. It is
also clear that $dC^0_l(\fm_2,\fm_2)$ is one-dimensional and generated by
$de_l\,=\,[e_l,-]$ for $l\geq 1$. Observe that $[e_1,-]$ is zero on $e_1$ and 
non-trivial on all other $e_i$, that $[e_2,-]$ is zero on $e_2$, equal to a
constant $a$ on all $e_i$ with $i\geq 3$ and equal to $-a$ on $e_1$, while
$[e_i,-]$ for $i\geq 3$ is non-zero on $e_1$ and $e_2$ and zero on all others.

One sees that $Z^1_1(\fm_2,\fm_2)=dC^0_1(\fm_2,\fm_2)$. We therefore conclude 
that  

\begin{theo}
$${\rm dim}\,H^1_l(\fm_2,\fm_2)\,=\,\left\{\begin{array}{ccccc} 0 & {\rm if} &
l=1 & {\rm or} & l\leq -1 \\
1 & {\rm if} & l=0 & {\rm or} & l\geq 2 \end{array}\right.$$ 
\end{theo}

This theorem has been found independently by Dimitri Millionschikov in
\cite{Million}. 
\smallskip

In order to compute the bracket structure, we need explicit non-trivial 
cocycles. Observe that the (non zero) coboundary for $l\geq 3$ is given by 
$a_1\not=0$ and $a_2=a_1$. The explicit non-trivial cocycles are therefore:
\begin{itemize}
\item $l=0$: the coefficients are growing linearly $a:=a_1$, $a_2=2a$, 
$a_3=3a$ etc.
\item $l=2$: $b:=a_2\not= 0$ and $a_j=b$ for all $j\geq 3$.
\item $l\geq 3$: $a_1=:-\frac{c_l}{2}$ and $a_2=\frac{c_l}{2}$. Then $a_3=c_l$,
$a_4=c_l$, etc. 
\end{itemize}
We express the previous description by introducing generators:
\begin{itemize}
\item $l=0$: $\omega(e_k)=ke_k$ for all $k\geq 1$ (we took $a=1$).
\item $l=2$: $$\alpha(e_k)=\left\{\begin{array}{ccc} be_{k+2} & {\rm if} 
& k\geq 2 \\ 0 & {\rm if} & k=1 \end{array}\right.$$ 
\item $l\geq 3$: $$\gamma_l(e_k)=\left\{\begin{array}{ccc} c_le_{k+l} & 
{\rm if} & k\geq 3 \\  -\frac{c_l}{2}e_{l+1} & {\rm if} & k=1  \\
\frac{c_l}{2}e_{l+2} & {\rm if} & k=2 \end{array}\right.$$
\end{itemize}

It is well known that $H^*(\fg,\fg)$ carries a graded Lie algebra structure for
any Lie algebra $\fg$, and that $H^1(\fg,\fg)$ forms a graded Lie subalgebra.
Let us compute this bracket structure on our generators:

Given $a\in C^p(\fg,\fg)$ and $b\in C^q(\fg,\fg)$, define
$$ab(x_1,\ldots,x_{p+q-1})\,=\,\sum_{\sigma\in{\rm Sh}_{p,q}}(-1)^{{\rm sgn}\,
\sigma}a(b(x_{i_1},\ldots,x_{i_q}),x_{j_1}\ldots,x_{j_{p-1}})$$
for $x_1,\ldots,x_{p+q-1}\in\fg$. The bracket is then defined by
$$[a,b]\,=\,ab-(-1)^{(p-1)(q-1)}ba.$$
It thus reads on $H^1(\fg,\fg)$ simply
$$[a,b](x)\,=\,a(b(x))-b(a(x)).$$

We compute
\begin{eqnarray*}
\omega(\alpha(e_k))-\alpha(\omega(e_k))&=&\left\{\begin{array}{ccc}
0 & {\rm if} & k=1 \\ \omega(be_{k+2}) & {\rm if} & k\geq 2 \end{array}\right\}
- \alpha(ke_k)    \\
&=&\left\{\begin{array}{ccc}
0 & {\rm if} & k=1 \\ b(k+2)e_{k+2}-bke_{k+2} & {\rm if} & k\geq 2 
\end{array}\right.   \\
&=&\left\{\begin{array}{ccc}
0 & {\rm if} & k=1 \\ 2be_{k+2} & {\rm if} & k\geq 2 \end{array}\right.\,=\,
2\alpha(e_k).
\end{eqnarray*}

\begin{eqnarray*}
\omega(\gamma_l(e_k))-\gamma_l(\omega(e_k))&=&\left\{\begin{array}{ccc}
\omega(c_le_{k+l}) & {\rm if} & k=\geq 3 \\ 
\omega(-\frac{c_l}{2}e_{l+1}) & {\rm if} & k=1  \\
\omega(\frac{c_l}{2}e_{l+2}) & {\rm if} & k=2 \end{array}\right\}
- \left\{\begin{array}{ccc}
kc_le_{k+l}) & {\rm if} & k=\geq 3 \\ 
-k\frac{c_l}{2}e_{l+1} & {\rm if} & k=1  \\
k\frac{c_l}{2}e_{l+2} & {\rm if} & k=2 \end{array}\right\}   \\
&=&\left\{\begin{array}{ccc}
c_l(k+l)e_{k+l}-kc_le_{k+l} & {\rm if} & k=\geq 3 \\ 
\left(-\frac{c_l}{2}(l+1)+\frac{c_l}{2}\right)e_{l+1} & {\rm if} & k=1  \\
\left(\frac{c_l}{2}(l+2)-c_l\right)e_{l+2} & {\rm if} & k=2 \end{array}\right\}
\,=\,l\gamma_l(e_k).
\end{eqnarray*}

\begin{eqnarray*}
\alpha(\gamma_l(e_k))-\gamma_l(\alpha(e_k))&=&\left\{\begin{array}{ccc}
\alpha(c_le_{k+l}) & {\rm if} & k=\geq 3 \\ 
\alpha(-\frac{c_l}{2}e_{l+1}) & {\rm if} & k=1  \\
\alpha(\frac{c_l}{2}e_{l+2}) & {\rm if} & k=2 \end{array}\right\}
- \left\{\begin{array}{ccc}
\gamma_l(be_{k+2}) & {\rm if} & k=\geq 3 \\ 
0 & {\rm if} & k=1  \\
\gamma_l(be_{4}) & {\rm if} & k=2 \end{array}\right\}   \\
&=&\left\{\begin{array}{ccc}
bc_le_{k+l+2}-bc_le_{k+l+2} & {\rm if} & k=\geq 3 \\ 
-\frac{c_l}{2}be_{l+1+2}-0 & {\rm if} & k=1  \\
\frac{c_l}{2}be_{l+4}-bc_le_{l+4} & {\rm if} & k=2 \end{array}\right\}  \\
&=&\left\{\begin{array}{ccc}
0 & {\rm if} & k=\geq 3 \\ 
-\frac{c_l}{2}be_{l+1+2} & {\rm if} & k=1  \\
-\frac{c_l}{2}be_{l+4} & {\rm if} & k=2 \end{array}\right\} \\
&=&\left\{\begin{array}{ccc}
0 & {\rm if} & k=\geq 3 \\ 
-\frac{c_l}{2}be_{k+l+2} & {\rm if} & k=1,2. \end{array}\right.  \\
\end{eqnarray*}
This last cocycle is a coboundary, more precisely,
$$\alpha(\gamma_l(e_k))-\gamma_l(\alpha(e_k))\,=\,\left(\frac{c_l}{2}b\right)
[e_{l+2},-].$$
We conclude  
$$\alpha(\gamma_l(e_k))-\gamma_l(\alpha(e_k))\,=\,0$$
as cohomology classes. One easily computes that $\gamma_l$ and $\gamma_m$ 
commute. Therefore the bracket structure on $H^1(\fm_2,\fm_2)$ is described 
as follows:

\begin{theo}
$H^1(\fm_2,\fm_2)$ is a graded Lie algebra, generated in positive degrees by
$\omega$ (degree $0$), $\alpha$ (degree $2$) and $\gamma_l$ (degree $l\geq 3$)
such that $\omega$ acts as a grading operator on the trivial Lie algebra 
generated by $\alpha$ and the $\gamma_l$ for $l\geq 3$.
\end{theo} 

\section{The space $H^2(\fm_2,\fm_2)$}

\subsection{Cocycle identities}

For a $2$-cochain $\omega$, the cocycle identity reads
\begin{eqnarray*} 
\omega([e_i,e_j],e_k)+\omega([e_j,e_k],e_i)+\omega([e_k,e_i],e_j)\\
-[e_i,\omega(e_j,e_k)]-[e_j,\omega(e_k,e_i)]-[e_k,\omega(e_i,e_j)]=0.
\end{eqnarray*}
In the sequel, we will suppose $\omega$ homogeneous of weight $l\in\Z$ with 
$\omega(e_i,e_j)=a_{i,j}e_{i+j+l}$ for all $i,j\geq 1$. From the cocycle 
identity, we get the following equations on the coefficients $a_{i,j}$:
\begin{itemize}
\item[(a)] Setting $i=1$ and $j,k\geq 3$, we get for $j+k+l\geq 2$
$$(a_{j+1,k} + a_{j,k+1})e_{j+k+l+1}=(a_{j,k}-a_{k,1}\delta_{k+l,0}
-a_{k,1}\delta_{k+l,1}-
a_{1,j}\delta_{j+l,0}-a_{1,j}\delta_{j+l,1})e_{j+k+l+1},$$
and for $j+k+l=0,1$ (while there is no equation for $j+k+l<0$)
$$(a_{j+1,k} + a_{j,k+1})e_{j+k+l+1}=0.$$
\item[(b)] Setting $i=1$, $j=2$, and $k\geq 3$, we get for $k+l\geq 2$, 
$$(a_{3,k}+a_{k+2,1}+a_{2,k+1})e_{k+l+3}=(a_{2,k}+a_{k,1}
-a_{1,2}\delta_{2+l,0}-a_{1,2}\delta_{2+l,1})e_{k+l+3},$$
while for $k+l=0$, we get
$$(a_{3,k}+a_{k+2,1}+a_{2,k+1})e_{3}=(a_{2,k}-a_{k,1}
-a_{1,2}\delta_{2+l,0}-a_{1,2}\delta_{2+l,1})e_{3},$$
and for $k+l=1$, we get
$$(a_{3,k}+a_{k+2,1}+a_{2,k+1})e_4=(a_{2,k}-a_{1,2}\delta_{2+l,0}-
a_{1,2}\delta_{2+l,1})e_4,$$
and for $k+l=-1,-2$, we get
$$(a_{3,k}+a_{k+2,1}+a_{2,k+1})e_{k+l+3}=(-a_{1,2}\delta_{2+l,0}
-a_{1,2}\delta_{2+l,1})e_{k+l+3}.$$
\item[(c)] If $i=2$, and $j,k\geq 3$, we get for $j+k+l\geq 3$ 
$$(a_{j+2,k}+a_{j,k+2})e_{j+k+l+2}=(a_{j,k}-a_{k,2}\delta_{k+l+1,0}
-a_{k,2}\delta_{k+l,0}-
a_{2,j}\delta_{j+l+1,0}-a_{2,j}\delta_{j+l,0})e_{j+k+l+2},$$
for $j+k+l=1$
$$(a_{j+2,k}+a_{j,k+2})e_3=(-a_{j,k}-a_{k,2}\delta_{k+l+1,0}
-a_{k,2}\delta_{k+l,0}-
a_{2,j}\delta_{j+l+1,0}-a_{2,j}\delta_{j+l,0})e_{3},$$
for $j+k+l\leq -2$, there is no equation, and for $j+k+l=0,-1,2$, we have
$$(a_{j+2,k}+a_{j,k+2})e_{j+k+l+2}=(-a_{k,2}\delta_{k+l+1,0}
-a_{k,2}\delta_{k+l,0}-a_{2,j}\delta_{j+l+1,0}-
a_{2,j}\delta_{j+l,0})e_{j+k+l+2}.$$
\item[(d)] If $i,j,k\geq 3$, we get 
$$0=(-a_{j,k}\delta_{j+k+l,1}-a_{j,k}\delta_{j+k+l,2}-a_{k,i}\delta_{i+k+l,1}
-a_{k,i}\delta_{i+k+l,2}-a_{i,j}\delta_{i+j+l,1}-a_{i,j}\delta_{i+j+l,2})
e_{i+j+k+l}.$$
\end{itemize}
In equation (d), at most two terms can be non-zero for a given $l$ as $i$, 
$j$ and $k$ must be pairwise distinct.


Let us now compute the $2$-coboundaries: a cocycle $\omega$ is a coboundary
in case there exists a $1$-cochain $\alpha$ such that
$$\omega(e_i,e_j)\,=\,\alpha([e_i,e_j])-[e_i,\alpha(e_j)]+[e_j,\alpha(e_i)].$$
As $\omega$ is homogeneous of weight $l$, $\alpha$ will be, and we set 
$\alpha(e_i)=a_ie_{i+l}$ for all $i\geq 1$. Then the previous equation gives:
\begin{itemize}
\item[(e)] Suppose $i=1$ and $j\geq 3$, then
$$a_{1,j}e_{j+l+1}\,=\,(a_{j+1}-a_j-a_1\delta_{l,0}-a_1\delta_{l+1,0})
e_{j+l+1}.$$
This equation makes sense only if $j+l\geq 2$. For $j+l=0,1$, one obtains
$$a_{1,j}\,=\,a_{j+1}.$$
\item[(f)] Suppose $i=1$ and $j=2$, then for $l\geq 2$
$$a_{1,2}\,=\,a_3-a_2+a_1,$$
while for $l=-1,-2$, one gets $a_{1,2}\,=\,a_3$, for $l=0$, one gets 
$a_{1,2}\,=\,a_3-a_2-a_1$, and for $l=1$, one gets 
$a_{1,2}\,=\,a_3-a_2$.
\item[(g)] Suppose $i=2$ and $j\geq 3$, then for $j+l\geq 2$
$$a_{2,j}\,=\,a_{j+2}-a_j(1-\delta_{j+l,2})+a_j\delta_{j+l,1}
-a_2(\delta_{l,0}+\delta_{l,-1}),$$
while for $j+l=0,-1$, one gets $a_{2,j}\,=\,a_{j+2}$, and for $j+l=1$, one gets
$$a_{2,j}\,=\,a_{j+2}+a_j.$$
\item[(h)] For $i,j\geq 3$ with $i+j+l\geq 1$, $i+l\geq 1$ and $j+l\geq 1$, 
the coboundary equation reads
$$a_{i,j}\,=\,a_j\,(\delta_{j+l,1}+\delta_{j+l,2})-a_i(\delta_{i+l,1}+
\delta_{i+l,2}).$$ 
\end{itemize}

Now stably, i.e. for a fixed $l$ and $j,k>>0$, we have just the following 
system of equations:
\begin{itemize}
\item[($\alpha$)] $a_{3,k}+a_{k+2,1}+a_{2,k+1}\,=\,a_{2,k}+a_{k,1}-a_{1,2}
\delta_{l,-2}-a_{1,2}\delta_{l,-1}$
\item[($\beta$)] $a_{j+1,k}+a_{j,k+1}\,=\,a_{j,k}$
\item[($\gamma$)] $a_{j+2,k}+a_{j,k+2}\,=\,a_{j,k}$
\end{itemize}

Equation ($\alpha$) means that the $1$- and $2$-coefficients determine the 
$3$-coefficients. Equation ($\beta$) implies that the differences of adjacent 
$3$- (resp. $4$-) coefficients determine the $4$- (resp. $5$-)
coefficients. But equation ($\gamma$) implies that differences of next to 
adjacent $3$-coefficients determine the $5$-coefficents directly. We get
stably on the one hand 
$$a_{5,k}\,=\,a_{4,k}-a_{4,k+1}\,=\,(a_{3,k}-a_{3,k+1})-(a_{3,k+1}-a_{3,k+2})
\,=\,a_{3,k}-2\,a_{3,k+1}+a_{3,k+2},$$
and on the other hand
$$a_{5,k}\,=\,a_{3,k}-a_{3,k+2},$$
thus we conclude that for $l$ big enough $a_{3,k+1}=a_{3,k+2}$. 
Even if $k>>0$, we take $j=3$ in order to get these equations, thus there
are extra terms (coming from equations (c)) for $j=-l$ and $j=-l-1$, i.e. 
in case $l=-3$ and $l=-4$. In all other weights, we will finally (i.e.
for $k>>0$) have the conclusion $a_{3,k+1}=a_{3,k+2}$.

But now when the $3$-coefficients are stably equal, 
the $4$-coefficients are stably $0$, and so are
all higher coefficients. This limits considerably the choice of possible
cocycles, at least stably. For example, let us suppose $l\geq -2$. In this 
case, equations (e) and (f) show that we can add coboundaries in order to have
all $1$-coefficients equal to zero. It is clear from equations (e), (f), (g) 
and (h) that once the $1$-coefficients are set to zero, the $2$-coefficients 
and higher coefficients cannot be changed by addition of a coboundary, because
this would change the $1$-coefficients, too.

(a), (b) and (c) then show that we have
the system of equations 
\begin{itemize}
\item[($\alpha'$)] $a_{3,k}+a_{2,k+1}\,=\,a_{2,k}$
\item[($\beta$)] $a_{j+1,k}+a_{j,k+1}\,=\,a_{j,k}$
\item[($\gamma$)] $a_{j+2,k}+a_{j,k+2}\,=\,a_{j,k}$
\end{itemize}
for all $j,k\geq 3$. 
The system tells us that cocycles must have all $3$-coefficients equal,
all higher coefficients zero. Observe that the equations which determined 
the solutions for $\fm_0$ are a subset of the equations which must be satisfied
for $\fm_2$. We conclude that in weight $l\geq -2$, there are at most two
non-trivial families of true deformations: the $2$-family and the $3$-family.
Whether they give indeed rise to true deformations will be determined in 
later subsections by studying their Massey powers. 

\subsection{Coboundaries}

In this subsection, we show that the $2$-family is a coboundary in weights
$l=-1,0,1$. 

Indeed, in weight $l=-1$, the coboundary equations read
\begin{eqnarray*}
0&=&a_{1,j}=a_{j+1}-a_j-a_1\,\,\,\,\,\,\,\,\,j\geq 3 \\
0&=&a_{1,2}=a_3\\
1&=&a_{2,j}=a_{2+j}-a_j-a_2\,\,\,\,\,\,\,\,\,j\geq 4 \\
1&=&a_{2,3}=a_5-a_2
\end{eqnarray*}

Therefore, we conclude that the choice $a_3=0$, $a_4=a_1$, $a_5=2a_1$, 
$a_6=3a_1$, etc, with $a_1=a_2=1$, shows that the $2$-family is a coboundary 
in weight $l=-1$.

Now in weight $l=0$, the cobounadry equations read  
\begin{eqnarray*}
0&=&a_{i,j}=a_{j+1}-a_j-a_1\,\,\,\,\,\,\,\,\,j\geq 3 \\
0&=&a_{1,2}=a_3-a_2-a_1\\
1&=&a_{2,j}=a_{2+j}-a_j-a_2\,\,\,\,\,\,\,\,\,j\geq 3 
\end{eqnarray*}

Therefore, the choice $a_1=0$, $a_2=-1=a_3=a_4=\ldots$ shows that the 
$2$-family is a coboundary in weight $l=0$.

Finally, in weight $l=1$, the coboundary equations read:
\begin{eqnarray*}
0&=&a_{1,j}=a_{j+1}-a_j-a_1\,\,\,\,\,\,\,\,\,j\geq 3 \\
0&=&a_{1,2}=a_3-a_2\\
1&=&a_{2,j}=a_{2+j}-a_j\,\,\,\,\,\,\,\,\,j\geq 3 
\end{eqnarray*}

Here the choice $a_1=\frac{1}{2}$, $a_2=a_3$, $a_4=a_3+\frac{1}{2}$, 
$a_5=a_4+\frac{1}{2}$, etc shows that the $2$-family is a coboundary.

One easily sees that it is not a coboundary in all other weights $\geq -2$ 
(for this it is enough to check $l=2,-2$, because the coboundary equations 
stabilize for $l\geq 2$, and in these two cases, writing the $2$-family as a 
coboundary leads to a contradiction).

Let us summarize the discussion of sections $3.1$, $3.2$ and $3.4$ in the 
following theorem ($H^2_l({\mathfrak m}_2,{\mathfrak m}_2)=\{0\}$ for 
$l\leq -5$ will be shown in $3.4$):

\begin{theo}
$$\dim\,H^2_l({\mathfrak m}_2,{\mathfrak m}_2)=\left\{\begin{array}{ccc} 
0 & {\rm for} & l\leq -5 \\
1 &  {\rm for} & l=-1,0,1 \\ 
2 &  {\rm for} & l=-4,-3,-2\,\,\,\,\,{\rm or}\,\,\,\,\,l\geq 2 
\end{array}\right.$$
\end{theo}

In particular, $H^2({\mathfrak m}_2,{\mathfrak m}_2)$ is infinite 
dimensional, but $H^2_l({\mathfrak m}_2,{\mathfrak m}_2)$ is finite dimensional
for each fixed $l\in\Z$. This means that Fialowski-Fuchs' construction of a 
miniversal deformation \cite{FiaFuc} does work (cf section $7.4$ in 
\cite{FiaFuc}). In order to get hold of it, one would need to have some 
information on $H^3_l({\mathfrak m}_2,{\mathfrak m}_2)$, like for example the 
informations displayed in \cite{Million}. 

Theorem $4$ has been found independently by Dimitri 
Millionschikov in \cite{Million}.   

\subsection{Massey powers}

Observe that the Massey square does not involve
the bracket of the Lie algebra, so we get for $\fm_2$ the same Massey square as
for $\fm_0$. For example, the $2$-family has zero Massey square (as a cochain)
in all weights (but observe that the $2$-family is not necessarily a cocycle in
all weights). We will examine the $3$-family in positive or zero weight in the 
following proposition.

An important point is that for $\fm_0$, we had restrictions on the true 
deformations coming from the nullity of the Massey squares and higher Massey
powers. For $\fm_2$ here, we have more possibilities to compensate non-zero
Massey powers, so there are less restrictions. Most of the restrictions for
deformations of $\fm_2$ come already from the cocycle equations.

\begin{prop}
Let $\omega\in Z^2_l(\fm_2,\fm_2)$ be the homogeneous $2$-cocycle of weight 
$l\geq 0$ given by the $3$-family and
representing an infinitesimal deformation of $\fm_2$. Then $\omega$ can be 
prolongated to a formal deformation of $\fm_2$, i.e. all Massey powers 
$[\omega]^n\in H^3(\fm_2,\fm_2)$ of $\omega$ are trivial.
\end{prop}

\pr Recall that the homogeneous $2$-cocycle $\omega$ of weight $l$ is given by 
coefficients $a_{i,j}$ such that $\omega(e_i,e_j)=a_{i,j}e_{i+j+l}$. 
$\omega$ represents the $3$-family, thus $a_{i,j}\not= 0$ (up to antisymmetry) 
only for $i=2$ and $j\geq 5$ and $i=3$ and $j\geq 4$. The Massey square
of $\omega$ reads
$$M_{ijk}\,=\,a_{i,j}a_{i+j+l,k}\,+\,a_{j,k}a_{j+k+l,i}\,+\,
a_{k,i}a_{k+i+l,j}.$$
We will always suppose $i<j<k$, up to anti-symmetry. Using $a_{i,j}\not= 0$ 
(up to antisymmetry) only for $(i=2,j\geq 5)$ and $(i=3,j\geq 4)$, we
obtain as only possibly non-zero Massey squares $M_{2jk}$, $j,k\geq 4$, and 
$M_{3jk}$, $j,k\geq 4$. The squares $M_{3jk}$, $j,k\geq 4$ are zero because
of the restriction $l\geq 0$; indeed,
$$M_{3jk}\,=\,a_{3,j}a_{3+j+l,k}\,+\,a_{k,3}a_{k+3+l,j}\,=\,a_{3+j+l,k}\,+\,
a_{k+3+l,j},$$
and $l\geq 0$, $j,k\geq 4$ imply that $a_{3+j+l,k}=a_{k+3+l,j}=0$.

The squares $M_{2jk}$, $j,k\geq 4$ are zero for $j\geq 4$, because then   
$$M_{2jk}\,=\,a_{2,j}a_{2+j+l,k}\,+\,a_{k,2}a_{k+2+l,j},$$
and once again, $l\geq 0$, $j,k\geq 4$ imply that $a_{2+j+l,k}=a_{k+2+l,j}=0$.

Therefore, the only Massey squares we have to compensate are $M_{23k}$, 
$k\geq 4$. We then introduce a homogeneous $2$-cochain $\alpha$ of weight $2l$
with $\alpha(e_i,e_j)=b_{i,j}e_{i+j+2l}$. We have for $l\geq -1$
$$d\alpha(e_2,e_j,e_k)\,=\,(b_{2+j,k}-b_{k+2,j}-b_{j,k})e_{j+k+2l+2},$$
meaning $d\alpha(e_2,e_3,e_k)\,=\,(b_{5,k}-b_{k+2,3}-b_{3,k})e_{k+2l+5}$.
We may then compensate the Massey square by just the $3$-column of 
$b$-coefficients. This ensures that at most the $2$- and $3$-columns for the 
$a$- and the $b$-coefficients are non-zero. 

Now suppose by induction that we have already compensated all Massey powers 
up to some level in such a way that at most the $2$- and $3$-columns for the 
coefficients of the intervening cochains are non-zero. 
Then we go on to compute the next Massey power
$$N_{ijk}\,=\,\beta(\gamma(e_i,e_j),e_k)\,+\,\gamma(\beta(e_i,e_j),e_k)\,+\,
{\rm cycl.},$$
where ``cycl.'' means cyclic permutations in $i,j,k$ and $\beta$ and $\gamma$
are some $2$-cochains satisfying the above restrictions. The weights of the 
cochains $\beta$ and $\gamma$ are positive or zero. Thus by compensating one 
step further, we will reproduce
cochains such that at most the $2$- and $3$-columns for the coefficients are 
non-zero. This ends the inductive step.\fin

Let us summerize what we said about true deformations in weight $l\geq 0$:

\begin{prop}
In weight $l\geq 0$, the only non-trivial cocycles are given by (linear 
combinations of) the $2$- and the $3$-family, but the $2$-family is a 
coboundary in weights $l=0,1$. The $2$-family gives rise to 
a true deformations (its Massey square is zero as a cochain), while the 
$3$-family gives rise to a formal deformation.
\end{prop}

We will be more specific about the convergence of this formal deformation and
about the $\N$-graded Lie algebras to which $\fm_2$ deforms in weight $0$ in a
later subsection.  

\subsection{Cocycles in weight $l\leq -5$}

Let us show in this section that there are no non-trivial $2$-cocycles in 
weight $l\leq -5$. This is somewhat surprising; we interprete it as being the 
fact that the cocycle equations for $\fm_2$ are very restrictive.

First of all, equations (e) mean that we can compensate the coefficients 
$a_{1,j}$ for $j+l\geq 0$ by a suitable coboundary. Observe that $a_{1,j}$ 
does not make sense for $j+l\leq -1$ as $a_{1,j}$ is the coefficient in front 
of $e_{j+1+l}$, so it can be set to zero. Therefore we will suppose 
in the following that $a_{1,j}=0$ for all $j\geq 2$. Thus, by antisymmetry, all
coefficients involving an index $1$ are zero.    

With this in mind, the cocycle equations (a) and (b) become more simple:

\begin{itemize}
\item  $a_{3,k}+a_{2,k+1}\,=\,a_{2,k}$
\item $a_{j+1,k}+a_{j,k+1}\,=\,a_{j,k}$  
\end{itemize}

for $k\geq 3$, $k+l\geq 2$, resp. $j+k+l\geq 2$, $j,k\geq 3$.  

Let us write down the cocycle equations of type (c) with $j=3$ 
(this is the case of interest for the reasoning which eliminates higher 
non-zero terms) and $k\geq 4$:

\begin{eqnarray*}
-l-4\leq k\leq -l-3: a_{5,k}&=& -a_{3,k+2} \\
k=-l-2: a_{5,k}&=& -a_{3,k}-a_{3,k+2} \\
k=-l-1: a_{5,k}&=& -a_{3,k+2}-a_{k,2} \\
k=-l: a_{5,k}&=& a_{3,k}-a_{3,k+2}+a_{2,k} \\
k\geq -l+1: a_{5,k}&=& a_{3,k}-a_{3,k+2} 
\end{eqnarray*}

Thus, for $k\geq -l+1$, we have on the one hand 
$a_{5,k}= a_{3,k}-a_{3,k+2}$, and on the other hand (for $k\geq -l-1$)
$$a_{5,k}\,=\,a_{4,k}-a_{4,k+1}\,=\,(a_{3,k}-a_{3,k+1})-(a_{3,k+1}-a_{3,k+2})
\,=\,a_{3,k}-2\,a_{3,k+1}+a_{3,k+2},$$
and one deduces $a_{3,k+1}=a_{3,k+2}$ for all $k\geq -l+1$. We call this 
coefficient $x:=a_{3,k+1}=a_{3,k+2}$ for all $k\geq -l+1$.

The equation $a_{5,-l}\,=\,a_{3,-l}-a_{3,-l+2}+a_{2,-l}$ and the equation
$a_{5,-l}\,=\,a_{3,-l}-2a_{3,-l+1}+a_{3,-l+2}$ imply that 
$2a_{3,-l+2}\,=\,2a_{3,-l+1}+a_{2,-l}$, and therefore with $a:=a_{2,-l}$,
we get $x\,=\,a_{3,-l+1}+\frac{a}{2}$. 

\noindent{\bf Step 1:}\quad Using these equations, we fill in the table of 
coefficients $a_{i,j}$ starting from high $k$ values:

\vspace{1cm}
\hspace{3cm}
\begin{tabular}{|c|c|c|c|c|c} \hline
       & $2$      &  $3$   &  $4$   &  $5$   \\  \hline\hline
$-l$   &  $a$   &      &      &       \\  \hline
$-l+1$ &        &  $x-\frac{a}{2}$&$-\frac{a}{2}$&$-\frac{a}{2}$    \\  \hline
$-l+2$ &      &  $x$   &  $0$   &  $0$   \\  \hline
$-l+3$ &      &  $x$   &  $0$   &  $0$   \\  \hline
$-l+4$ &      &  $x$   &  $0$   &  $0$   \\  \hline
$-l+5$ &      &  $x$   &  $0$   &  $0$   \\  \hline
\end{tabular}
\vspace{1cm}

The $-\frac{a}{2}$ will repeat itself to the right of the table, meaning
$a_{4+r,-l+1}=-\frac{a}{2}$ for all $r$. But $a_{-l+1,-l+1}=0$ by antisymmetry,
thus $a=0$. 

\noindent{\bf Step 2:}\quad When we call $a_{3,-l}=:y$, the new table looks 
like: 

\vspace{1cm}
\hspace{3cm}
\begin{tabular}{|c|c|c|c|c|c} \hline
     & $2$      &  $3$     &  $4$    &  $5$   \\  \hline\hline
$-l$   &  $0$   &  $y$   &  $y-x$&  $y-x$  \\  \hline
$-l+1$ &        &  $x$   &  $0$  &  $0$    \\  \hline
$-l+2$ &        &  $x$   &  $0$  &  $0$   \\  \hline
$-l+3$ &        &  $x$   &  $0$  &  $0$   \\  \hline
$-l+4$ &        &  $x$   &  $0$  &  $0$   \\  \hline
$-l+5$ &        &  $x$   &  $0$  &  $0$   \\  \hline
\end{tabular}
\vspace{1cm}

Once again, continuing the line with $y-x$ to the right,
when we hit the diagonal, we get $y=x$. 

\noindent{\bf Step 3:}\quad When we call $a_{3,-l-1}=:a$, the new table looks 
like:

\vspace{1cm}
\hspace{3cm}
\begin{tabular}{|c|c|c|c|c|c} \hline
     & $2$      &  $3$     &  $4$    &  $5$   \\  \hline\hline
$-l-1$ &  $a$   &  $a-x$ &  $a-x$& $a-x$   \\  \hline
$-l$   &  $0$   &  $x$   &  $0$  &  $0$  \\  \hline
$-l+1$ &  $-x$  &  $x$   &  $0$  &  $0$    \\  \hline
$-l+2$ &  $-2x$ &  $x$   &  $0$  &  $0$   \\  \hline
$-l+3$ &  $-3x$ &  $x$   &  $0$  &  $0$   \\  \hline
$-l+4$ &  $-4x$ &  $x$   &  $0$  &  $0$   \\  \hline
$-l+5$ &  $-5x$ &  $x$   &  $0$  &  $0$   \\  \hline
\end{tabular}
\vspace{1cm}

The same argument as before gives us here $x=a$. 

\noindent{\bf Step 4:}\quad This time, call $a_{3,-l-2}=:y$, then we get 
by the equation $a_{5,-l-2}\,=\,a_{3,-l-2}-a_{3,-l}$ that $a_{5,-l-2}\,=\,-y-x$
and $a_{5,-l-2}\,=\,-x+y\,=\,a_{4,-l-2}\,=\,-a_{3,-l-1}+a_{3,-l-2}$. One 
concludes $y=0$.

\noindent{\bf Step 5:}\quad Now write the new table:

\vspace{1cm}
\hspace{3cm}
\begin{tabular}{|c|c|c|c|c|c} \hline
     & $2$      &  $3$     &  $4$    &  $5$   \\  \hline\hline
$-l-2$ &  $x$   &  $0$   &  $-x$ &  $-x$   \\  \hline
$-l-1$ &  $x$   &  $x$   &  $0$  &  $0$   \\  \hline
$-l$   &  $0$   &  $x$   &  $0$  &  $0$  \\  \hline
$-l+1$ &  $-x$  &  $x$   &  $0$  &  $0$    \\  \hline
$-l+2$ &  $-2x$ &  $x$   &  $0$  &  $0$   \\  \hline
$-l+3$ &  $-3x$ &  $x$   &  $0$  &  $0$   \\  \hline
$-l+4$ &  $-4x$ &  $x$   &  $0$  &  $0$   \\  \hline
$-l+5$ &  $-5x$ &  $x$   &  $0$  &  $0$   \\  \hline
\end{tabular}
\vspace{1cm}

Finally, hitting once again the diagonal shows that $x=0$. In order to conclude
that all coefficients must be zero, it suffices to show that $a_{4,-l-3}=0$.
This follows from the (a) equation (with $j=3$, $k=-l-3$): 
$a_{4,-l-3}=-a_{3,-l-2}=0$. $a_{4,-l-3}=0$ suffices, because $a_{i,j}$ can only
be non-zero starting from $i+j+l\geq 0$, i.e. $i=2$ and $j\geq -l-1$, $i=3$ and
$j\geq -l-2$, $i=4$ and $j\geq -l-2$ and so on. 

We summarize in the following

\begin{prop}
There are no non-trivial $2$-cocycles in weight $l\leq -5$.
\end{prop}

\subsection{True deformations in weights $l=-1$ and $l=-2$}

Again, by the same reasoning as before, all coefficients involving an index 
$1$ can be set to zero (up to addition of coboundaries). 

The (a) and (b) equations are like in the general case. The (c) equations are
not yet modified (only for $l=-3$ and $l=-4$). There is no non trivial (d) 
equation yet. 

We are thus still in the range of validity of the reasoning which shows that 
there are as only possibly non-trivial cocycles the $2$- and the $3$-family.

The $2$-family is still a cocycle of Massey square zero (as a cochain).
The only thing which may be different here is the proof that the $3$-family
gives still rise to a formal deformation.

The first steps are like in the proof of proposition $1$: 
the only Massey squares we have to compensate are $M_{23k}$, 
$k\geq 4$. We then introduce a homogeneous $2$-cochain $\alpha$ of weight $2l$
with $\alpha(e_i,e_j)=b_{i,j}e_{i+j+2l}$. We have for $l\geq -2$, 
$j,k\geq 3$, $j<k$:
$$d\alpha(e_2,e_3,e_k)\,=\,b_{5,k}-b_{k+2,3}-b_{3,k}+
\delta_{k+2l,0}b_{k,2}+\delta_{5+2l,1}b_{2,3}+\delta_{5+2l,2}b_{2,3}.$$
As for the $3$-family $b_{2,3}=0$, this reads more simply:
$$d\alpha(e_2,e_3,e_k)\,=\,b_{5,k}-b_{k+2,3}-b_{3,k}+
\delta_{k+2l,0}b_{k,2}.$$
We may choose to compensate once again just by the $3$-column, i.e. we may
set $b_{5,k}=b_{k,2}=0$ for all $k$. This ensures that at most the $2$- and 
$3$-columns for the $a$- and the $b$-coefficients are non-zero.

The next Massey power is then the Massey cube:  

\begin{eqnarray*}
N_{ijk}&=&\alpha(\omega(e_i,e_j),e_k)\,+\,\omega(\alpha(e_i,e_j),e_k)\,+\,
{\rm cycl.}  \\
&=&a_{i,j}b_{i+j+l,k}\,+\,b_{i,j}a_{i+j+2l,k}\,+\,{\rm cycl.}.
\end{eqnarray*} 





We see that the terms we have to compensate are once again of type $N_{23k}$ 
(up to antisymmetry). We will have more and more Massey 
powers to compensate. This can be achieved by a growing, but finite number of 
cochains at each level. On the other hand, this process will not stop.
We therefore get:

\begin{prop}
In weight $l=-1,-2$, the only homogeneous $2$-cocycles are the $2$- and the 
$3$-family. The $2$-family is a coboundary in weight $l=-1$.
The $2$-family is of Massey square zero (as a cochain), and gives 
thus rise to a true deformation in weight $l=-2$. 
The $3$-family has zero Massey powers, and 
gives rise to a formal deformation with non-zero contributions at each level.
\end{prop}

\subsection{True deformations in weights $l=-3$ and $l=-4$}

Let us write down the cocycle equations. The important equations are those of 
type (c). They read:
$$a_{j+2,k}+a_{j,k+2}\,=\,a_{j,k} -a_{k,2}\delta_{k+l+1,0}
-a_{k,2}\delta_{k+l,0}-
a_{2,j}\delta_{j+l+1,0}-a_{2,j}\delta_{j+l,0}.$$

In weight $l=-3$, this means for $j=3$ and $k\geq 4$ that
$$a_{5,k}\,=\,a_{3,k}-a_{3,k+2}-a_{2,3}.$$
Compare this equation to 
$$a_{5,k}\,=\,a_{3,k}-2a_{3,k+1}+a_{3,k+2},$$
which follows as usually from the (a) equations. In conclusion, we get:
$$-2a_{3,k+2}\,=\,a_{2,3}-2a_{3,k+1}.$$
This means once again that the differences of $3$-coefficients are constant, 
and thus that the $4$-coefficients are equal, while the $5$-coefficients are 
zero. More precisely
$$2a_{4,k+1}\,=\,2(a_{3,k+1}-a_{3,k+2})\,=\,a_{2,3},$$
and therefore $a_{4,k+1}=\frac{a_{2,3}}{2}$. Either $a_{4,k+1}\not=0$ and
we get a family with non-zero coefficients in the first three columns, or 
$a_{4,k+1}=0$, i.e. $a_{2,3}=0$, and we get the $3$-family.  

Observe that the $2$-family does not satisfy the cocycle identities in weight
$l\leq -3$. Indeed, for $j,k\geq 3$    
$$a_{j+2,k}+a_{j,k+2}\,=\,a_{j,k} -a_{k,2}\delta_{k+l+1,0}
-a_{k,2}\delta_{k+l,0}-a_{2,j}\delta_{j+l+1,0}-a_{2,j}\delta_{j+l,0},$$
and for $k>>0$, all terms are zero, but one of the form $a_{2,j}$. This is a 
contradiction.

It remains thus (a linear combination of) the $3$- and the $4$-family. The
$3$-family is of Massey square zero in weight $l=-3$ (see the $\fm_0$-case !).

Let us turn to weight $l=-4$. Once again we look at a $2$-cocycle $\omega$
given by coefficients $a_{i,j}$ such that $a_{1,k}=0$ for all $k\geq 2$, which
we can achieve possibly by adding a coboundary, cf equations (e). We 
cannot exploit independently equations (f) and (g), because in these equations
the same coefficients occur. 

Let us write down low degree (a) equations: 
$$a_{j+1,k}+a_{j,k+1}\,=\,a_{j,k},$$
for $j,k\geq 3$. We therefore have for example $a_{3,4}=a_{3,5}$. 
The (b) equations read 
\begin{itemize}
\item $k=3$: $a_{2,4}=0$.
\item $k=4$: $a_{3,4}+a_{2,5}=a_{2,4}=0$. 
\item $k=5$: $a_{3,5}+a_{2,6}=a_{2,5}$.
\item $k\geq 6$: $a_{3,k}+a_{2,k+1}=a_{2,k}$.
\end{itemize}
And the (c) equations, which are the most interesting, read for $j=3$:
\begin{itemize}
\item $k=4$: $a_{5,4}+a_{3,6}=a_{3,4}-a_{4,2}-a_{2,3}=a_{3,4}-a_{2,3}$.
\item $k=5$: $a_{3,7}=a_{3,5}-a_{2,3}$. 
\item $k\geq 6$: $a_{5,k}+a_{3,k+2}=a_{3,k}-a_{2,3}$.
\end{itemize} 
The (d) equations are still void. 

Let us now start a table with the coefficients $a_{i,j}$ which varify these 
equations. First of all, we call $a:=a_{2,3}$, and $a_{3,4}=:b$. Then 
on the one hand $-a_{4,5}=b-a-a_{3,6}$, and on the other hand 
$a_{4,5}=b-a_{3,6}$. This gives $a_{3,6}=b-\frac{a}{2}$, and 
$a_{4,5}=\frac{a}{2}$. 

Now let us perform the same trick as in the other cases: on the one hand, we
have $a_{5,k}=a_{3,k}-a_{3,k+2}-a_{2,3}$, and on the other hand, we have 
$a_{5,k}=a_{3,k}-2a_{3,k+1}+a_{3,k+2}$ by the (a) equations, for $k\geq 6$.
We get thus $a_{3,k+1}-a_{3,k+2}=\frac{a_{2,3}}{2}$, i.e. the differences
of the $3$-coefficients, which determine the $4$-coefficients, are constant,
and therefore the $5$-coefficients zero. We now display the table:

\vspace{1cm}
\hspace{3cm}
\begin{tabular}{|c|c|c|c|c|c} \hline
    &   2    &  3     &  4    &  5   \\  \hline\hline
$3$ &  $a$   &        &       &       \\  \hline
$4$ &  $0$   &  $b$   &       &        \\  \hline
$5$ &  $-b$   &  $b$   &$\frac{a}{2}$&         \\  \hline
$6$ &  $-2b$  &  $b-\frac{a}{2}$   &$\frac{a}{2}$&  $0$    \\  \hline
$7$ &$-3b+\frac{a}{2}$&  $b-\frac{2a}{2}$   &$\frac{a}{2}$&  $0$   \\  \hline
$8$ &$-4b+\frac{3a}{2}$&  $b-\frac{3a}{2}$   &$\frac{a}{2}$&  $0$   \\  \hline
$9$ &$-5b+\frac{6a}{2}$&  $b-\frac{4a}{2}$   &$\frac{a}{2}$&  $0$   \\  \hline
$10$ &$-6b+\frac{10a}{2}$& $b-\frac{5a}{2}$  &$\frac{a}{2}$&  $0$   \\  \hline
\end{tabular}
\vspace{1cm}

We see that a $2$-parameter family is building up. The remaining question is
whether the Massey powers are zero, i.e. whether the family gives rise to a
true or formal deformation. We will consider the two cases $a=0$ and $b=0$
separately. For $b=0$, we have (a multiple of) the $4$-family (up to a
non-zero coefficient $a_{2,3}$). One easily
varifies that the additional non-zero coefficient $a_{2,3}$ does not change 
the Massey square zero character of the $4$-family in weight $l=-4$ (cf the
$\fm_0$-case). For $a=0$, we have the $3$-family which has non-zero Massey
squares. We compute that $M_{234}=0$, $M_{235}=0$, but $M_{23k}\not=0$ for
$k\geq 6$, that $M_{245}\not=0$, but $M_{24k}=0$ for $k\geq 6$, that 
$M_{25k}\not=0$ for $k\geq 6$, that $M_{26k}=0$ for $k\geq 7$, that
$M_{34k}\not=0$ for $k\geq 5$ and finally that $M_{35k}=0$ for $k\geq 6$.
These are all ordered Massey squares which are possibly non-zero.

We have thus a finite family of non-zero Massey squares which can be 
compensated by a finite sum of coboundaries. These give then rise to a finite
number of higher dimensional Massey powers, which can also be compensated in
the usual way. All in all we get a formal deformation.

\begin{prop}
In weights $l=-3$ and $l=-4$, the $3$-family and the $4$-family (and their 
linear combinations) are the only $2$-cocycles. In weight $l=-3$, the 
$3$-family gives a true and the $4$-family a formal deformations, whereas in
weight $l=-4$, the $4$-family gives a true and the $3$-family a formal 
deformation.
\end{prop}

\subsection{Identification of the deformations in weight $l=0$}   

We have seen in one of the previous sections that there is exactly one 
non-trivial cocycle in weight $l=0$. It is given by the $3$-family. 
We then examined Massey powers, and found that the $3$-family has Massey 
squares at each step and gives finally rise to a formal deformation. 
Let us identify in this section the Lie algebras to which $\fm_2$ deforms.

Consider the deformation given by the $3$-family. The corresponding 
deformation $\fm_2^2(t)$ reads (up to antisymmetry):
$$[e_1,e_j]_t = e_{j+1}\,\,\,\,\,\,\forall j\geq 2,$$    
$$[e_2,e_j]_t = e_{j+2}+t(1-(j-4))e_{j+2}\,\,\,\,\,\,\forall j\geq 4,$$
$$[e_2,e_j]_t = te_{j+3}\,\,\,\,\,\,\forall j\geq 4.$$
We already saw that this deformation has Massey corrections in any power of
$t$, so that it is a formal deformation. Let us show that it gives a 
non-converging deformation. Indeed, if it were converging, the limiting
object would be an $\N$-graded Lie algebra with one-dimensional graded 
components, generated in degrees $1$ and $2$. But by the classification 
theorem (Theorem p. 2 in \cite{Fia1}), $\fm_2^2(t)$ must be isomorphic to
$L_1$. This is obviously not the case, as $\fm_2^2(t)$ has a codimension
$3$ abelian ideal, whereas $L_1$ does not have any abelian ideal.

Therefore we arrive at the conclusion:

\begin{prop}
The deformations of $\fm_2$ in weight $l=0$ described in the following way:
the only non-trivial $2$-cocycle leads to
a formal non-converging deformation. In particular, $\fm_2$ does not deform
to any other $\N$-graded Lie algebra with one-dimensional graded 
components, generated in degrees $1$ and $2$. In particular, it does not deform
to $L_1$.
\end{prop}

\end{document}